\documentclass{article}
\usepackage[margin=1in]{geometry}
\usepackage{amsmath,amssymb,amsthm}
\usepackage{mathtools}
\usepackage{array}
\usepackage{bm}
\usepackage{enumitem}
\usepackage[page]{appendix}
\usepackage{natbib}
\usepackage[hidelinks]{hyperref}
\usepackage[nameinlink]{cleveref}

\newtheorem{theorem}{Theorem}[section]
\newtheorem{lemma}{Lemma}[section]
\newtheorem*{corollary*}{Corollary}
\newtheorem*{remark*}{Remark}

\DeclareMathOperator{\rank}{rank}

\DeclareMathOperator{\col}{col}

\renewcommand{\Pr}{\mathbb{P}}

\newcommand{\RR}{\mathbb{R}}

\newcommand{\set}[1]{\mathcal{#1}}
\newcommand{\covars}{[d]}

\newcommand{\projmat}{\Pi}
\newcommand{\ortprojmat}{\Pi^{\perp}}

\newcommand{\mineig}{\lambda_{\min}}
\newcommand{\betamin}{\beta_{\min}}

\crefname{appendix}{appendix}{appendices}
\Crefname{appendix}{Appendix}{Appendices}

\title{Best Subset Selection in Linear Regression:\\ Fixed-Design Error Bounds and Insights for Random Designs}

\author{
    Maxim Fedotov\\[0.75em]
    Department of Economics and Business\\
    Universitat Pompeu Fabra\\
    Barcelona, Spain\\
    \texttt{maxim.fedotov@upf.edu}
}

\date{July 2026}

\begin{document}

\maketitle

\begin{abstract}
We study exact support recovery by best subset selection in linear regression under a fixed design. For a known support size $s$ and ambient dimension $d$, we derive a non-asymptotic upper bound on the probability that best subset selection fails to recover the true support. The bound is expressed through a deterministic subset-separability parameter, which measures how well the true support can be distinguished from competing supports after projection. The result holds for all sample sizes $n$ exceeding a certain sufficient threshold which we state explicitly in terms of the signal-to-noise ratio, the subset separability of the realized design, and a logarithmic factor of order $\ln s + \ln(d - s)$. In contrast to random-design analyses, no full log-combinatorial term over the candidate support class appears. We discuss how such terms may reappear when the design is random and the separability parameter must be controlled uniformly over many competing subsets. The fixed-design formulation and the proof strategy also indicate settings in which the effective complexity of best subset selection may be reduced, for instance, under structured designs or restricted candidate subset classes.
\end{abstract}

\section{Introduction} 
\label{sec:intro}

We revisit the classical best subset selection problem in linear regression. Suppose that the $n$-dimensional outcome vector $Y$ satisfies
\begin{equation}\label{eq:linear-regression}
    Y = X \beta + W,
\end{equation}
where $X \in \RR^{n\times d}$ is a fixed design matrix, $\beta \in \RR^d$ is an unknown vector of coefficients, and $W \sim \mathcal{N}(0, \sigma^2 I_n)$ is a noise vector. We assume that the model is sparse in that $\beta$ contains zeros.
The main goal of best subset selection is to recover $\set{S} \coloneqq \{j \colon \beta_j \ne 0\}$ based on observations $(X, Y)$. For a method selecting some subset $\widehat{\set{S}}$, the error probability $\Pr\{\widehat{\set{S}}\ne\set{S}\}$ is a standard measure of performance. Its dependence on the sample size $n$, ambient dimension $d$, sparsity level $s \coloneqq |\set{S}|$, signal strength $\betamin \coloneqq \min_{j\in\set{S}} |\beta_j|$, and noise level $\sigma^2$ is a central question in sparse regression.

We suppose that the true support size $s \le \min\{n,d\}$ is known and consider the best subset selector
\begin{equation}\label{eq:best-subset}
    \widehat{\set{S}}_{n, s} 
    \in
    \underset{\set{M} \subseteq [d] \colon |\set{M}| = s}{\arg\min} \ 
    \min_{\beta'_{\set{M}} \in \RR^{|\set{M}|}} \ 
    \|Y - X_{\set{M}} \beta'_{\set{M}}\|_2^2,
\end{equation}
where $X_{\set{M}}$ denotes the submatrix of $X$ containing the columns indexed by $\set{M}$.
Assuming known $s$ provides necessary conditions for support recovery in unknown-sparsity settings, and thus is common and widely studied in the literature; see, for example, \citet{wainwright2009information, gamarnik2022sparse,shen2013constrained}.

Under a Gaussian random design, where the rows of $X$ are drawn i.i.d. from $\mathcal{N}(0,\Sigma)$, the seminal work of \cite{wainwright2009information} establishes that, for consistent recovery of a fixed true support $\set{S}$ by $\widehat{\set{S}}_{n,s}$, it suffices that
\begin{equation}\label{eq:sample-size-gaussian-ensemble}
    n 
    \gtrsim 
    \max\left\{
        \ln\binom{d - s}{s},\ 
        \frac{\sigma^2\ln (d - s)}{\rho(\Sigma) \cdot \betamin^2}
    \right\},
\end{equation}
where $\rho(\Sigma) \coloneqq \min_{\set{C} \subseteq [d]\colon |\set{C}| = s,\ \set{C} \neq \set{S}} \mineig(\Sigma_{\set{S} \setminus \set{C}, \set{S} \setminus \set{C}} - \Sigma_{\set{S} \setminus \set{C}, \set{C}} \Sigma_{\set{C}, \set{C}}^{-1}\Sigma_{\set{C}, \set{S} \setminus \set{C}})$, with $\mineig(\,\cdot\,)$ denoting the minimum eigenvalue of a symmetric matrix, is a subset-separability parameter. The former term in \eqref{eq:sample-size-gaussian-ensemble} is log-combinatorial: under $s\le d/2$, $\ln\binom{d-s}{s}$ is the logarithm of the number of competing supports obtained by replacing all $s$ active covariates by inactive ones. 
Since $\ln\binom{d-s}{s}\le s\ln(d-s)$, the signal-dependent term dominates whenever $\rho(\Sigma) \cdot \betamin^2/\sigma^2\lesssim 1/s$. In stronger-signal or better-separated regimes, especially when $s$ grows, the log-combinatorial term becomes the visible sample-size cost.

Related information-theoretic results for sparse linear models and noisy compressed sensing were established by \citet{akccakaya2009shannon}, \citet{fletcher2009necessary}, \citet{reeves2008sampling}, and \citet{wang2010information}. In this literature, necessary and sufficient sample size conditions for reliable support recovery typically feature either log-combinatorial terms involving the ambient dimension $d$, as in \eqref{eq:sample-size-gaussian-ensemble}, or related information-theoretic quantities reflecting the size of the model class. Similar quantities also appear in more recent work concerning computational and statistical phase transitions in sparse regression \citep{bertsimas2020sparse,gamarnik2022sparse}.

Several fixed-design analyses of $\ell_0$-based variable selection are particularly close to the present note. 
\citet{shen2013constrained} obtain finite-sample exponential error bounds for $\ell_0$-constrained least squares in terms of a degree-of-separation quantity measuring the residual prediction gap between the true model and its closest competitors. 
\citet{guo2020best} develop a related fixed-design theory for best subset selection based on an identifiability margin, and \citet{roy2025understanding} refine this viewpoint by showing that the required margin is also affected by geometric complexities of residualized signals and spurious projections. 
A key difference is that the main separation or margin quantities in these works combine the realized geometry of $X$ with the signal $\beta$ as they are built from residualized linear predictions.
In contrast, the subset-separability parameter used here is a signal-free empirical quantity, and the signal strength enters separately through $\betamin^2$.

Another related result is due to \citet{rognonvael2026improving}, who establish an error probability bound for $\ell_0$-penalized regression with external information under unknown sparsity. From the perspective of this note, a conceptually similar empirical subset-separability quantity appears in their analysis, and their bound leads to comparable sufficient sample size scaling. However, the $\ell_0$-penalty is essential for their bound and is subject to certain assumptions, which are not present in the current note as our focus is the known-sparsity case.

The present note is complementary in emphasis and formulation. For the fixed-cardinality best subset selector, we derive a non-asymptotic error probability bound and the associated sufficient sample size condition in a way that separates realized-design geometry, signal strength, noise level, and support combinatorics. This allows us to draw novel insights into the statistical complexity of best subset selection under a random design.

More specifically, we first establish a pairwise error probability bound comparing the true support with an arbitrary competing subset of the same size. 
The bound depends on the residualized design geometry through the unexplained variance of the missed active covariates and the corresponding part of the signal. 
This gives a natural definition of uniform subset separability, $r_{\set{S},n}$, which we use to obtain an explicit error probability bound for $\widehat{\set{S}}_{n,s}$ by combining the pairwise bounds.
The resulting sufficient sample size condition depends on $\betamin^2/\sigma^2$, $r_{\set{S},n}$, and a logarithmic factor of order $\ln s+\ln(d-s)$. 
Thus, no log-combinatorial term such as $\ln\binom{d-s}{s}$ in \eqref{eq:sample-size-gaussian-ensemble} or $\ln\binom{d}{s}$ appears explicitly in the fixed-design sufficient condition.
The result also applies beyond the usual sparse regime $s\le d/2$, including the non-sparse case $s>d/2$.

Finally, we discuss how our fixed-design result can be used in random-design settings and what implications it yields. 
We present two sources of reduced effective complexity for best subset selection. 
The first is geometric: the design may reduce the number of effectively distinct projection directions generated by competing subsets.
We illustrate this by a latent-factor design in which inactive covariates lie in a common low-dimensional factor subspace, so that the random-design cost of controlling subset separability is governed by the factor dimension and the number of missed active variables, rather than by all possible false replacements. 
The second is combinatorial and arises when the candidate support class is restricted due to additional information or structural modeling constraints such as grouped variables. 
In both cases, the fixed-design viewpoint makes explicit that the relevant complexity enters through the empirical subset-separability parameter of the realized design.

The rest of the note is organized as follows. \Cref{sec:setup} introduces the fixed-design model and the quantities used to obtain the theoretical results, including the pairwise comparison statistic and the subset-separability parameter. \Cref{sec:results} contains the pairwise and global error probability bounds for best subset selection. \Cref{sec:discussion} discusses the implications of the fixed-design result for random designs, structured covariate distributions, and restricted candidate classes.

\subsection*{Notation}

We let $[d] \coloneqq \{1, \ldots, d\}$. For a subset $\set{M} \subseteq \covars$, we write $\beta_\set{M}$ for the corresponding subvector of $\beta$. The complement of $\set{M}$ in $\covars$ is denoted by $\set{M}^c$, and $|\set{M}|$ denotes its cardinality. For two nonnegative functions $f$ and $g$, we write $f \gtrsim g$ if there exists a universal constant $c > 0$ such that $f \ge c\, g$. For a matrix $M$, $\rank(M)$ denotes its rank, $\col(M)$ is its column span, and $\mineig(M)$ is the smallest eigenvalue if $M$ is symmetric. The Moore--Penrose pseudoinverse is denoted by $(\,\cdot\,)^+$. In the fixed-design results, $\Pr$ denotes probability with respect to the noise vector $W$.

\section{Setup}
\label{sec:setup}

Best subset selection can be analyzed through the prism of pairwise comparisons between the true support and the competing subsets. For the support $\set{S}$ and all competitors $\set{C} \ne \set{S}$, we define the corresponding pairwise statistic as
\begin{equation}\label{eq:Delta-hat}
    \widehat{\Delta}_{\set{C}, \set{S}} 
    \coloneqq 
    \|Y - X_\set{C} \hat{\beta}_\set{C}\|_2^2 
    - \|Y - X_\set{S} \hat{\beta}_\set{S}\|_2^2
\end{equation}
where $\hat{\beta}_\set{M} \in \arg\min_{\beta'_\set{M} \in \RR^{|\set{M}|}} \ \|Y - X_\set{M} \beta'_\set{M}\|_2^2$. Simple algebra shows that, for any $\set{C}$,
\begin{equation}
    \widehat{\Delta}_{\set{C}, \set{S}} = \|\ortprojmat_\set{C} Y\|_2^2 - \|\ortprojmat_\set{S} Y\|_2^2
\end{equation}
where $\ortprojmat_\set{M} \coloneqq I_n - X_\set{M} (X_\set{M}^T X_\set{M})^+ X_\set{M}^T$ is the orthogonal projection operator onto the complement to the column space of $X_\set{M}$.

By the definition in \eqref{eq:best-subset}, the best subset selector $\widehat{\set{S}}_{n, s}$ fails to recover $\set{S}$ correctly and unambiguously if and only if there exists at least one subset $\set{C} \ne \set{S}$ such that $\widehat{\Delta}_{\set{C}, \set{S}} \le 0$. Hence, the associated error probability satisfies
\begin{equation}\label{eq:error-probability}
    \Pr\{\widehat{\set{S}}_{n, s} \neq \set{S}\} 
    \le  
    \Pr\left[
        \cup_{\set{C} \subseteq [d] \colon |\set{C}| = s, \ \set{C} \ne \set{S}}
        \{\widehat{\Delta}_{\set{C}, \set{S}} \le 0\}
    \right].
\end{equation}
In \Cref{sec:results}, we establish a non-asymptotic upper bound for the error probability for a pairwise comparison between the support $\set{S}$ and an arbitrary competing subset $\set{C} \ne \set{S}$ of size $s$. Then, based on this result, the overall error probability \eqref{eq:error-probability} is bounded via the union bound. 

The pairwise bound is based on a convenient decomposition for $\widehat{\Delta}_{\set{C}, \set{S}}$ featuring a pairwise ``signal-separation'' quantity 
$\frac{1}{n}\|\ortprojmat_\set{C} X_{\set{S} \setminus \set{C}} \beta_{\set{S} \setminus \set{C}}\|_2^2 = \|G_{\set{S},\set{C}}^{1/2} \beta_{\set{S} \setminus \set{C}}\|_2^2$,
where we define
\begin{equation}\label{eq:G}
    G_{\set{S}, \set{C}} 
    \coloneqq 
    \frac{1}{n}X_{\set{S} \setminus \set{C}}^T \ortprojmat_\set{C} X_{\set{S} \setminus \set{C}}.
\end{equation}
The matrix $G_{\set{S}, \set{C}}$ represents residual variance when regressing $X_{\set{S} \setminus \set{C}}$ onto $X_{\set{C}}$. Consider a hypothetical case where all columns of $X_{\set{S} \setminus \set{C}}$ lie in the column span of $X_{\set{C}}$. In this case, $\set{S}$ and $\set{C}$ are not distinguishable from the perspective of the linear model, and $G_{\set{S}, \set{C}} = 0$. As for another example, suppose that $X_{\set{S} \setminus \set{C}}^TX_\set{C} = 0$. Then, $G_{\set{S}, \set{C}} = \frac{1}{n} X_{\set{S} \setminus \set{C}}^TX_{\set{S} \setminus \set{C}}$. Thus, higher variance among variables in $\set{S} \setminus \set{C}$ yields larger magnitude of $G_{\set{S}, \set{C}}$.

The overall error probability bound is, in its nature, uniform over competing subsets and signals with the minimum magnitude given by $\betamin$. So, it requires uniform control on the pairwise signal separation which we address through a notion of subset separability. Similarly to the random-design subset separability utilized by \cite{wainwright2009information}, we define the fixed-design version as
\begin{equation}\label{eq:r-S}
    r_{\set{S}, n}
    \coloneqq
    \min_{\set{C} \subseteq \covars \colon |\set{C}| = s,\ \set{C} \ne \set{S}}
    \mineig(G_{\set{S}, \set{C}}).
\end{equation}
This parameter provides a useful lower bound for the pairwise signal separation,
\[
    \|G_{\set{S},\set{C}}^{1/2}
    \beta_{\set{S} \setminus \set{C}}\|_2^2
    \ge
    r_{\set{S}, n} \cdot \betamin^2 \cdot |\set{S} \setminus \set{C}|.
\]

\Cref{theorem:probability-of-error} assumes $r_{\set{S}, n} > 0$. Even though this condition is not needed for the pairwise bound, it is natural for the uniform guarantee. If $r_{\set{S}, n}=0$, then there exists at least one competing subset $\set{C}$ and infinitely many coefficient vectors $\beta \not \equiv 0$ such that $X_{\set{S} \setminus \set{C}}\beta_{\set{S} \setminus \set{C}}$ lies in the span of $X_\set{C}$, and thus $\|G_{\set{S},\set{C}}^{1/2}\beta_{\set{S} \setminus \set{C}}\|_2^2 = 0$. In this case, omitting the variables in $\set{S} \setminus \set{C}$ induces no deterministic loss in the criterion used by best subset selection. Other coefficient directions may still be separated, but this cannot be captured by the error probability bound which utilizes $r_{\set{S}, n}$, and thus is uniform across all coefficient directions.

We also assume that $X_\set{S}$ is full rank, that is, $\rank(X_\set{S}) = s$.
If this assumption is violated, then the effective support size is lower than $s$, which would affect some concentration inequalities used to obtain the non-asymptotic error bounds. In the classical example where the rows of $X$ arise from a multivariate Gaussian distribution, $X_\set{S}$ is full rank with probability one if the distribution of the corresponding covariates is marginally non-degenerate and $s \le \min\{n, d\}$.

\section{Fixed-design error probability bounds}
\label{sec:results}

The following lemma shows that the probability of $\set{S}$ not being preferred over $\set{C}$ by the best subset selection scheme from \eqref{eq:best-subset} can be controlled. For the proof and the discussion regarding the stated conditions, see \Cref{appendix:lemma:pairwise-statistic}.
\begin{lemma} \label{lemma:pairwise-statistic}
    Let $c = 4 \, (3 + \sqrt{5})$. If $\rank(X_\set{S}) = s$, then for all subsets $\set{C} \ne \set{S}$ of size $s$ satisfying $\|G_{\set{S}, \set{C}}^{1/2} \beta_{\set{S} \setminus \set{C}}\|_2^2 \ge \frac{c \sigma^2 \cdot |\set{S} \setminus \set{C}|}{n}$,
    \[\Pr\{\widehat{\Delta}_{\set{C}, \set{S}} \le 0\} \le 3 \exp\left\{- \frac{\|G_{\set{S}, \set{C}}^{1/2} \beta_{\set{S} \setminus \set{C}}\|_2^2}{c \sigma^2} \cdot n\right\}.\]
\end{lemma}

Note that $\|G_{\set{S}, \set{C}}^{1/2} \beta_{\set{S} \setminus \set{C}}\|_2^2 \ge \mineig(G_{\set{S}, \set{C}}) \cdot \|\beta_{\set{S} \setminus \set{C}}\|_2^2 \ge r_{\set{S}, n} \cdot \betamin^2 \cdot |\set{S} \setminus \set{C}|$ for all $\set{C}$,
where the last inequality follows from the definition of $r_{\set{S}, n}$ in \eqref{eq:r-S}.
This implies an upper bound which is expressed directly through the number of missing active covariates and the problem parameters of interest, that are the signal-to-noise ratio and subset separability.
\begin{corollary*}[\Cref{lemma:pairwise-statistic}]\hypertarget{corollary:lemma:pairwise-statistic}{}
    If $\rank(X_\set{S}) = s$ and $r_{\set{S}, n} \cdot n \ge c \frac{\sigma^2}{\betamin^2}$, then for all subsets $\set{C} \ne \set{S}$ of size $s$,
    \[
        \Pr\{\widehat{\Delta}_{\set{C}, \set{S}} \le 0\} 
        \le 
        3 \exp\left\{
            -\frac{1}{c} 
            \cdot \frac{\betamin^2}{\sigma^2} 
            \cdot |\set{S} \setminus \set{C}| 
            \cdot r_{\set{S}, n} 
            \cdot n
        \right\}.
    \]
\end{corollary*}
Replacing the pairwise signal separation by the uniform lower bound with respect to $\beta_{\set{S} \setminus \set{C}}$ directions makes the probability bound in \hyperlink{corollary:lemma:pairwise-statistic}{the corollary} less sharp and the resulting sufficient condition stronger than in \Cref{lemma:pairwise-statistic}. However, it allows us to arrive at the uniform error probability bound in \Cref{theorem:probability-of-error}. As it turns out, keeping $|\set{S} \setminus \set{C}|$ is still important as it yields a tighter final bound.

Based on the pairwise result in \hyperlink{corollary:lemma:pairwise-statistic}{the corollary}, the error probability of the best subset selection rule from \eqref{eq:best-subset} can be bounded above given that the sample size surpasses a certain threshold. Both the sufficient sample size condition and the error probability bound feature $\betamin \coloneqq \min_{j \in \set{S}} |\beta_j|$, so the result can be viewed as a uniform guarantee over signals whose minimum nonzero coefficient is at least of magnitude $\betamin$. This is a typical perspective in the literature on variable selection. 

The theorem below is the main fixed-design error bound of the note. Its role is to isolate the deterministic separability parameter through which the design enters the error probability. The implications of this formulation for random designs are discussed in \Cref{sec:discussion}. The case $s = d$ is trivial since it implies only one admissible support. So, all bounds involving $r_{\set{S},n}$ and $\ln(d - s)$ are understood in the nontrivial case $s < d$. The proof is presented in \Cref{appendix:theorem:probability-of-error}.

\begin{theorem}\label{theorem:probability-of-error}
    If $\rank(X_\set{S}) = s < d$, and $r_{\set{S}, n} > 0$, then there exists a universal constant $c$ such that for all $c_1 \in [0, 1)$ and all
    \begin{equation}\label{eq:sample-size}
        n 
        \ge 
        \frac{c}{1 - c_1} 
        \cdot \frac{\sigma^2}{\betamin^2} 
        \cdot \frac{1}{r_{\set{S}, n}} 
        \cdot \max\{
            \ln s + \ln(d - s) + \ln \min\{s, d - s\},\
            1
        \},
    \end{equation}
    the error probability of best subset selection is bounded as
    \begin{equation}\label{eq:probability-bound}
        \Pr\{\widehat{\set{S}}_{n, s} \ne \set{S}\} 
        \le 
        3 \exp\left\{
            -\frac{c_1}{c} 
            \cdot \frac{\betamin^2}{\sigma^2} 
            \cdot r_{\set{S}, n} 
            \cdot n
        \right\}.
    \end{equation}
\end{theorem}

\begin{remark*}
    Under $\rank(X_\set{S}) = s$, the condition $r_{\set{S}, n} > 0$ can be verified through a rank condition. For every competing subset $\set{C} \ne \set{S}$ of size $s$, suppose that
    \[
        \rank\left(X_{\set{S} \cup \set{C}}\right)
        =
        |\set{S} \setminus \set{C}| + \rank(X_\set{C}).
    \]
    This condition implies $r_{\set{S}, n} > 0$ in the fixed-design setup. If the rows of $X$ are drawn independently from an absolutely continuous distribution with respect to Lebesgue measure on $\RR^d$ and $n \ge s + \min\{s,\ d-s\}$, then the rank condition holds with probability one.
\end{remark*}

We next comment on several features of \Cref{theorem:probability-of-error}.

The proof of \Cref{theorem:probability-of-error} proceeds by applying the union bound to the pairwise error probabilities from \Cref{lemma:pairwise-statistic}. 
Nevertheless, no full log-combinatorial term such as $\ln\binom{d-s}{s}$ or $\ln\binom{d}{s}$ appears explicitly in either the sufficient condition \eqref{eq:sample-size} or the final probability bound \eqref{eq:probability-bound}. 
This is a consequence of the exponential dependence of the pairwise bound on $|\set{S}\setminus\set{C}|$, together with a careful control of the combinatorial summation over competing subsets. 
The relation between this fixed-design phenomenon and the log-combinatorial terms appearing in random-design analyses is discussed in \Cref{sec:discussion}.

The parameter $c_1$ reflects a tradeoff between the sufficient sample size condition and the error probability. Larger values of $c_1$ strengthen the exponential decay in \eqref{eq:probability-bound}, but also increase the required sample size in \eqref{eq:sample-size}. In turn, for large enough $n$, one may take $c_1$ arbitrarily close to one, tightening the bound. The signal-to-noise ratio $\betamin^2/\sigma^2$ and the subset separability $r_{\set{S}, n}$ enter both expressions in the expected direction: stronger signal and better separation make support recovery easier.

Finally, \Cref{theorem:probability-of-error} is not restricted to sparse regimes. It remains applicable for all $s \le \min\{n, d - 1\}$, including non-sparse settings such as $s > d/2$. In this case, the maximum amount of active covariates that a competing subset $\set{C}$ of size s can miss, denoted by $|\set{S} \setminus \set{C}|$, is $d - s$ rather than $s$. This is reflected by the factor $\min\{s, d - s\}$ in \eqref{eq:sample-size}.

\section{Implications for random designs}
\label{sec:discussion}

\subsection{From fixed-design bounds to random-design guarantees}

\Cref{theorem:probability-of-error} is a fixed-design, non-asymptotic statement: the parameters $n,d,s,\set{S},\betamin,\sigma$ and the design matrix $X$ are treated as given. Nevertheless, the result provides a useful way to interpret the statistical complexity of best subset selection. Apart from the signal-to-noise ratio $\betamin^2/\sigma^2$, the sufficient condition in \eqref{eq:sample-size} is governed by the realized subset-separability parameter $r_{\set{S},n}$. In particular, the condition requires
\[
    n
    \gtrsim
    \frac{\sigma^2}{\betamin^2}
    \cdot
    \frac{\ln s+\ln(d-s)}
         {r_{\set{S},n}}.
\]
Thus, for a given design matrix, exact recovery is controlled by how well the true support remains separated from competing supports after projection, rather than directly by the total number of competing subsets.

This perspective clarifies the relation with random-design results. In the Gaussian ensemble setting, necessary and sufficient conditions for reliable support recovery typically involve log-combinatorial terms such as $\ln\binom{d-s}{s}$, $\ln\binom{d}{s}$, or related information-theoretic quantities; see, for example, \citet{wainwright2009information,wang2010information,gamarnik2022sparse}. At first sight, this may appear to contrast with \eqref{eq:sample-size}, which contains no explicit log-combinatorial term. The distinction is that, under a random design, $r_{\set{S},n}$ is itself random. Hence, applying the fixed-design theorem to a random-design setting requires controlling the lower tail of this separability parameter, uniformly over the relevant competing subsets. This additional step is precisely where combinatorial complexity may reappear.

The passage from fixed-design to random-design guarantees can be formalized by conditioning on the realized design. Let $\Xi_n$ be a class of design matrices $X$ such that, on $\Xi_n$, the assumptions of \Cref{theorem:probability-of-error} hold and
\[
    r_{\set{S},n}(X) \ge \underline r_{\set{S},n},
\]
where $\underline r_{\set{S},n} > 0$ is deterministic. Given $\Xi_n$, the error probability can be bounded as follows:
\[
    \Pr\{\widehat{\set{S}}_{n, s}\ne \set{S}\}
    =
    \mathbb{E}_X\!\left[
        \Pr\{\widehat{\set{S}}_{n, s}\ne \set{S}\mid X\}
    \right]
    \le
    \sup_{X\in\Xi_n}
    \Pr\{\widehat{\set{S}}_{n, s}\ne \set{S}\mid X\}
    +
    \Pr\{X \in \Xi_n^c\}.
\]
Conditionally on $\Xi_n$, the fixed-design theorem can be applied with
$r_{\set{S},n}$ replaced by $\underline r_{\set{S},n}$. Hence, if the corresponding
sample size condition holds, we have
\[
    \Pr\{\widehat{\set{S}}_{n, s}\ne \set{S}\}
    \le
    3\exp\left\{
        -\frac{c_1}{c}
        \cdot \frac{\betamin^2}{\sigma^2}
        \cdot \underline r_{\set{S},n} n
    \right\}
    +
    \Pr\{X \in \Xi_n^c\}.
\]

Thus, extending the fixed-design result to a random-design setting reduces to proving a high-probability lower bound for the empirical subset-separability parameter. In applications, this lower-tail control usually imposes its own sample-size requirement. So, a sufficient condition for a random design is obtained by combining the fixed-design recovery condition, with $r_{\set{S},n}$ replaced by $\underline r_{\set{S},n}$, and the design-specific condition that makes $\Pr\{X \in \Xi_n^c\}$ small. The resulting sufficient sample size is then the maximum of these two requirements. Different random-design structures lead to different $\Xi_n$, lower bounds $\underline r_{\set{S},n}$, and probabilities $\Pr\{X \in \Xi_n^c\}$. 

The next two subsections illustrate two sources of reduced random-design complexity. The first is geometric: many false replacement variables may generate the same low-dimensional nuisance directions. The second is combinatorial: the candidate support class itself may be restricted.

\subsection{Latent-factor designs}

We illustrate the geometric source of reduced complexity by the stylized latent-factor model detailed in \Cref{appendix:latent-factor-example}.
For a fixed true support $\set{S}=\{1,\ldots,s\}$, suppose that the design matrix is generated as
\[
    X_\set{S}=ZA+\tau U,
    \qquad
    X_{\set{S}^c}=ZB,
\]
where $Z\in\RR^{n\times v}$ and $U\in\RR^{n\times s}$ have independent standard Gaussian entries, $A$ and $B$ are deterministic loading matrices, and $\tau>0$ is the idiosyncratic standard deviation.
Thus, active covariates contain both a common low-dimensional factor component and an idiosyncratic component, whereas all inactive covariates lie in the same factor subspace of rank at most $v$.

The role of the factor structure is as follows.
Consider a competing subset $\set{C} \neq \set{S}$.
Since $X_{\set{C}\setminus\set{S}}$ is contained in $\col(Z)$, the span of the competitor is contained in $\col(X_{\set{S}\cap\set{C}},Z)$.
Hence, for the purpose of lower-bounding the Gram matrix for the residualized $X_{\set{S}\setminus\set{C}}$, it is enough to project out the larger space $\col(X_{\set{S}\cap\set{C}},Z)$, which no longer depends on the particular false replacement set $\set{C}\setminus\set{S}$.
After this projection, the common factor component is removed and the remaining separation is provided by the idiosyncratic active variation $\tau U_{\set{S}\setminus\set{C}}$.

The calculation in \Cref{appendix:latent-factor-example} shows that the lower tail of the empirical subset-separability parameter is controlled at the scale of the idiosyncratic variance $\tau^2$. In particular, under proportionality conditions on $n,s$ and $v$ made explicit in the appendix, there exist constants $c_0, C_0>0$ such that
\[
    \Pr\left\{
        r_{\set{S}, n}
        \le
        c_0 \tau^2
    \right\}
    \le
    \exp\{-C_0 n\}.
\]
Combining this bound with the conditioning argument discussed before gives
\[
    \Pr\{\widehat{\set{S}}_{n, s}\ne \set{S}\}
    \le
    3\exp\left\{
        -\frac{c_1 c_0}{c}
        \cdot \frac{\betamin^2}{\sigma^2}
        \cdot \tau^2 \cdot n
    \right\}
    +
    \exp\{-C_0 n\}
\]
provided that both the conditions giving the lower-tail bound above and the fixed-design sample size condition with $\underline r_{\set{S},n}=c_0\tau^2$ hold.

Note that high correlation between active and inactive covariates through the common factor does not by itself destroy support recovery in this example, because missed active variables retain an idiosyncratic component that inactive variables cannot reproduce.
The size of this component determines the separation scale: if $\tau$ is small, the bound becomes correspondingly weak.

The source of savings is different. At distance $m$ from the true support, the usual count
$\binom{s}{m}\binom{d-s}{m}$ is replaced by $\binom{s}{m}$ in the lower-tail control of $r_{\set{S},n}$, because the false replacement variables affect the argument only through the common factor space of dimension at most $v$.
Consequently, the empirical subset-separability parameter can be controlled without paying the full log-combinatorial price associated with choosing false replacements among the $d-s$ inactive variables.

\subsection{Restricted candidate classes}

The second source of potential reduction of complexity is a restriction of the candidate model class. 
Suppose that, instead of searching over all subsets of size $s$, the admissible supports are restricted to a smaller family $\mathfrak{S}$ with $\set{S}\in\mathfrak{S}$. 
The corresponding separability parameter then is
\[
    r_{\set{S},n}(\mathfrak{S})
    \coloneqq
    \min_{\set{C}\in\mathfrak{S}\colon \set{C}\ne\set{S}}
    \mineig
    \left(
        \frac{1}{n}
        X_{\set{S}\setminus\set{C}}^T
        \ortprojmat_\set{C}
        X_{\set{S}\setminus\set{C}}
    \right).
\]
Such restrictions can arise from external information, prior knowledge, screening, or modeling constraints such as grouped variables. 
The fixed-design proof then ranges only over competitors in $\mathfrak{S}$: at distance $m$ from $\set{S}$, the unrestricted count $\binom{s}{m}\binom{d-s}{m}$ is replaced by the number of admissible competitors in $\mathfrak{S}$ at that distance. 
Under a random design, the same restriction also reduces the number of lower-tail events that must be controlled for empirical subset separability.

This effect is illustrated in \Cref{appendix:restricted-candidate-class} for a standard Gaussian design. 
There, if the restricted population-level separability $\rho_{\set{S}}(\mathfrak{S})$ is positive and $n \gtrsim \max\{s, \ln|\mathfrak{S}|\}$, then there exist constants $c,C>0$ such that
\[
    \Pr\left\{
        r_{\set{S},n}(\mathfrak{S})
        \le
        c\rho_{\set{S}}(\mathfrak{S})
    \right\}
    \le
    \exp\{-Cn\}.
\]
Combining this lower-tail bound with the fixed-design argument gives a random-design recovery guarantee once the fixed-design sample-size condition is also satisfied with $r_{\set{S},n}(\mathfrak{S})$ replaced by $c\rho_{\set{S}}(\mathfrak{S})$.

\section{Conclusion}
\label{sec:conclusion}

Overall, the fixed-design result should be interpreted as a way to separate the role of the realized design from the signal strength and other problem parameters in support recovery. In this formulation, the design enters the error probability through the subset-separability parameter, which measures how well the active variables remain separated from competing supports after residualization. This perspective suggests that sharper random-design guarantees may be possible in structured settings where subset separability can be controlled without paying the full combinatorial price. In the absence of such structural design features, the usual log-combinatorial terms may remain unavoidable.

\section*{Acknowledgements}

The author is grateful to David Rossell and Gábor Lugosi for their careful reading of the manuscript and helpful comments. This work was supported by the predoctoral program AGAUR-FISDUR (2023 FISDUR 00566) of the Secretariat of Universities and Research of the Department of Research and Universities of the Generalitat of Catalonia and the European Social Plus Fund.

\bibliography{ref}

\newpage

\begin{appendices}
\crefalias{section}{appendix}

The appendix contains additional notation, proofs of the main fixed-design results, two random-design illustrations showing how empirical subset separability can be controlled under additional structure, and auxiliary results.

\section{Additional notation}
\label{appendix:notation}

We collect here several pieces of notation used in the appendix. If $M_1, \ldots, M_k$ are matrices with the same number of rows, then $[M_1, \ldots, M_k]$ denotes their column-wise concatenation. The notation $\ortprojmat_{[M_1, \ldots, M_k]}$ denotes the orthogonal projector onto the orthogonal complement of $\col([M_1, \ldots, M_k])$. For symmetric matrices $A$ and $B$, we write $A \succeq B$ when $A - B$ is positive semidefinite.

For nonnegative sequences $a_n$ and $b_n$, the notation $a_n \asymp b_n$ means that there exist two universal constants $0 < c < C < \infty$ and $n_0 \in \mathbb{N}$ such that $c b_n \le a_n \le C b_n$ for all $n \ge n_0$. We write $a_n \ll b_n$ if $a_n / b_n \to 0$ when $n \to \infty$.

In the random-design illustrations in \Cref{appendix:latent-factor-example,appendix:restricted-candidate-class}, probabilities are taken with respect to the random design matrix.

\section{Proof of \texorpdfstring{\Cref{lemma:pairwise-statistic}}{}}
\label{appendix:lemma:pairwise-statistic}

The proof strategy is similar to that of the one of Theorem 1 in \cite{wainwright2009information}. The probability of interest is bounded above by the union bound concerning two events, and two intermediate results that control the corresponding probabilities are presented. \Cref{lemma:projected-error-difference} is almost a direct adjustment to Lemma 3 of \cite{wainwright2009information} and \Cref{lemma:competitor-projected-outcome-deviation} is inspired by Lemma 4 of that same paper. The novelty in our proof is that it applies to the fixed-design case and allows the submatrices of $X$ given by subsets of $s$ columns, except $\set{S}$, to be rank-deficient.

Note that $\|\ortprojmat_\set{S} Y\|_2^2 = \|\ortprojmat_\set{S} W\|_2^2$. So, we can rewrite the pairwise comparison statistic for a pair of subsets $(\set{C}, \set{S})$ as
\begin{equation}
    \widehat{\Delta}_{\set{C}, \set{S}} 
    = 
    \|\ortprojmat_\set{C} Y\|_2^2 - \|\ortprojmat_\set{S} W\|_2^2 
    = 
    \|\ortprojmat_\set{C} Y\|_2^2 - \|\ortprojmat_\set{C} W\|_2^2 + \|\ortprojmat_\set{C} W\|_2^2 - \|\ortprojmat_\set{S} W\|_2^2,
\end{equation}
where we added and subtracted the same quantity $\|\ortprojmat_\set{C} W\|_2^2$.
For all $\delta \ge 0$, we have 
\begin{equation}\label{eq:Delta-hat-prob}
    \Pr\{\widehat{\Delta}_{\set{C}, \set{S}} \le 0\} 
    \le
    \Pr\{\|\ortprojmat_\set{C} Y\|_2^2 - \|\ortprojmat_\set{C} W\|_2^2 \le \delta \cdot n\} 
    + \Pr\{\|\ortprojmat_\set{S} W\|_2^2 - \|\ortprojmat_\set{C} W\|_2^2 > \delta \cdot n\}.
\end{equation}
where the last inequality follows from the union bound.

The following lemma establishes an upper bound on the latter summand. The proof is given in \Cref{appendix:lemma:projected-error-difference}.
\begin{lemma}
\label{lemma:projected-error-difference}
    If $\rank(X_\set{S}) = s$, then for all subsets $\set{C}$ of size $s$ and all $\delta~\ge~\frac{8 \sigma^2 \cdot |\set{S} \setminus \set{C}|}{n}$,
    \[
        \Pr\{
            \|\ortprojmat_\set{S} W\|_2^2 - \|\ortprojmat_\set{C} W\|_2^2 
            > 
            \delta \cdot n
        \} 
        \le 
        2 \exp\left\{-\frac{\delta}{8 \sigma^2} \cdot n\right\}.
    \]
\end{lemma}

The following lemma controls the former probability in (\ref{eq:Delta-hat-prob}). Recall that we defined 
\[
    G_{\set{S}, \set{C}} 
    \coloneqq 
    \frac{1}{n} 
    X_{\set{S} \setminus \set{C}}^T 
    \ortprojmat_\set{C} 
    X_{\set{S} \setminus \set{C}}
\]
and let $G_{\set{S}, \set{C}}^{1/2}$ be a matrix square root of $G_{\set{S}, \set{C}}$, which always exists as  $G_{\set{S}, \set{C}}$ is symmetric. The proof is presented in \Cref{appendix:lemma:competitor-projected-outcome-deviation}
\begin{lemma} \label{lemma:competitor-projected-outcome-deviation}
    For all subsets $\set{C}$ and all $c_1 > 0$, there exists $\delta_{c_1} \coloneqq \frac{1}{1 + c_1} \cdot \|G_{\set{S}, \set{C}}^{1/2} \beta_{\set{S} \setminus \set{C}}\|_2^2$ such that
    \[
        \Pr\left\{
            \|\ortprojmat_\set{C} Y\|_2^2 - \|\ortprojmat_\set{C} W\|_2^2 
            \le 
            \delta_{c_1} \cdot n
        \right\} 
        \le 
        \exp\left\{-\frac{c_1^2 \cdot \delta_{c_1}}{8 \sigma^2 \cdot (1 + c_1)} \cdot n \right\}.
    \]
\end{lemma}

Finally, we use \Cref{lemma:projected-error-difference} and \Cref{lemma:competitor-projected-outcome-deviation} together to bound the right hand side of \eqref{eq:Delta-hat-prob} above by choosing a particular value for $\delta$ which satisfies the requirements of the lemmas. To match the constants in both bounds, we take $c_1 = \frac{1 + \sqrt{5}}{2}$ in \Cref{lemma:competitor-projected-outcome-deviation} which yields $\delta = \frac{1}{1 + c_1} \cdot \|G_{\set{S}, \set{C}}^{1/2} \beta_{\set{S} \setminus \set{C}}\|_2^2$. 
\Cref{lemma:projected-error-difference} requires $\delta^*~\ge~\frac{8 \sigma^2 \cdot |\set{S} \setminus \set{C}|}{n}$, and thus we impose the following condition:
\[
    \|G_{\set{S}, \set{C}}^{1/2} \beta_{\set{S} \setminus \set{C}}\|_2^2 
    \ge 
    \frac{c\sigma^2 \cdot |\set{S} \setminus \set{C}|}{n}.
\]
where we denoted $c \coloneqq 8 (1 + c_1)$.
Under this condition, adding up the bounds established in \Cref{lemma:projected-error-difference} and \Cref{lemma:competitor-projected-outcome-deviation} yields the result:
\begin{equation}
    \Pr\{\widehat{\Delta}_{\set{C}, \set{S}} \le 0\} 
    \le 
    3 \exp\left\{
        -\frac{\|G_{\set{S}, \set{C}}^{1/2} \beta_{\set{S} \setminus \set{C}}\|_2^2}{c \sigma^2} 
        \cdot n
    \right\}.
\end{equation}

\section{Proofs of auxiliary results for \texorpdfstring{\Cref{lemma:pairwise-statistic}}{}}

\subsection{Proof of \texorpdfstring{\Cref{lemma:projected-error-difference}}{}}
\label{appendix:lemma:projected-error-difference}

The proof follows the argument of Lemma 3 in \citet{wainwright2009information}, adapted to the present fixed-design setting and to allow rank-deficient competing submatrices. We kindly refer the reader to the original paper in case of doubt.
    
Applying the Pythagorean theorem and some projection identities yields
\[
    \|\ortprojmat_\set{S} W\|_2^2 - \|\ortprojmat_\set{C} W\|_2^2 
    = 
    \|(\projmat_\set{C} - \projmat_{\set{S} \cap \set{C}}) W\|_2^2 
    - \|(\projmat_\set{S} - \projmat_{\set{S} \cap \set{C}}) W\|_2^2, 
\]
where $\projmat_\set{M} \coloneqq X_\set{M} (X_\set{M}^T X_\set{M})^+ X_\set{M}^T$ is a projection matrix onto the column span of $X_\set{M}$. Therefore, denoting $Z_{\set{C}, \set{S}} \coloneqq \frac{1}{\sigma^2}\|(\projmat_\set{C} - \projmat_{\set{S} \cap \set{C}}) W\|_2^2$ and similarly $Z_{\set{S}, \set{C}}$, we can rewrite the probability of interest as
\[
    \Pr\{
        \|\ortprojmat_\set{S} W\|_2^2 - \|\ortprojmat_\set{C} W\|_2^2 
        > 
        \delta \cdot n
    \} 
    = 
    \Pr\left\{
        Z_{\set{C}, \set{S}} - Z_{\set{S}, \set{C}} 
        > 
        \frac{\delta}{\sigma^2} \cdot n
    \right\}.
\]
Let $R_{\set{C}, \set{S}} \coloneqq \rank(X_\set{C}) - \rank(X_{\set{S} \cap \set{C}})$, and $R_{\set{S}, \set{C}}$ is defined similarly. $Z_{\set{C}, \set{S}}$ and $Z_{\set{S}, \set{C}}$ are central chi-square random variables with degrees of freedom $R_{\set{C}, \set{S}}$ and $R_{\set{S}, \set{C}}$ respectively. Denote $k_{\set{C}, \set{S}} \coloneqq R_{\set{C}, \set{S}} - R_{\set{S}, \set{C}} = \rank(X_\set{C}) - \rank(X_\set{S})$. Then,
\[
    \Pr\left\{
        Z_{\set{C}, \set{S}} - Z_{\set{S}, \set{C}} 
        > 
        \frac{\delta}{\sigma^2} \cdot n
    \right\} 
    = 
    \Pr\left\{
        Z_{\set{C}, \set{S}} - R_{\set{C}, \set{S}} 
        - (Z_{\set{S}, \set{C}} - R_{\set{S}, \set{C}}) 
        > 
        \frac{\delta}{\sigma^2} \cdot n - k_{\set{C}, \set{S}}
    \right\}.
\]
Since $\rank(X_\set{S}) = s$ by assumption, we have $k_{\set{C}, \set{S}} \le 0$. Hence,
\[
    \Pr\left\{
        Z_{\set{C}, \set{S}} - R_{\set{C}, \set{S}} 
        - (Z_{\set{S}, \set{C}} - R_{\set{S}, \set{C}}) 
        > 
        \frac{\delta}{\sigma^2} \cdot n - k_{\set{C}, \set{S}}\right
    \} 
    \le 
    \Pr\left\{
        Z_{\set{C}, \set{S}} - R_{\set{C}, \set{S}} 
        - (Z_{\set{S}, \set{C}} - R_{\set{S}, \set{C}}) 
        > 
        \frac{\delta}{\sigma^2} \cdot n
    \right\}.
\]
So, by the union bound,
\[
    \Pr\left\{
        Z_{\set{C}, \set{S}} - R_{\set{C}, \set{S}} 
        - (Z_{\set{S}, \set{C}} - R_{\set{S}, \set{C}}) 
        > 
        \frac{\delta}{\sigma^2} \cdot n
    \right\} 
    \le 
    \Pr\left\{
        Z_{\set{C}, \set{S}} - R_{\set{C}, \set{S}} 
        \ge 
        \frac{\delta}{2\sigma^2} \cdot n
    \right\} 
    + \Pr\left\{
        Z_{\set{S}, \set{C}} - R_{\set{S}, \set{C}} 
        \le 
        -\frac{\delta}{2\sigma^2} \cdot n
    \right\}
\]

The former term can be bounded by \Cref{lemma:chi-square-normalized-upper-tail-bound}. That is, if $\delta \ge \frac{8\sigma^2 \cdot R_{\set{C}, \set{S}}}{n}$, we have 
\[
    \Pr\left\{
        Z_{\set{C}, \set{S}} - R_{\set{C}, \set{S}} 
        \ge 
        \frac{\delta}{2\sigma^2} \cdot n
    \right\} 
    \le 
    \exp\left\{-\frac{\delta}{8\sigma^2} \cdot n\right\}
\]
The requirement on $\delta$ is satisfied since $R_{\set{C}, \set{S}} \le |\set{C} \setminus \set{S}| = |\set{S} \setminus \set{C}|$ and $\delta~\ge~\frac{8 \sigma^2 \cdot |\set{S} \setminus \set{C}|}{n}$ by assumption.

The latter term admits an even tighter bound by \Cref{lemma:chi-square-normalized-lower-tail-bound}. Given that $\rank(X_\set{S}) = s$ is assumed, we have $R_{\set{S}, \set{C}} > 0$ for all $\set{C} \neq \set{S}$. Applying
\Cref{lemma:chi-square-normalized-lower-tail-bound} with $t = \frac{\delta^2}{16\sigma^4 \cdot R_{\set{S},\set{C}}^2} \cdot n^2$
gives
\[
    \Pr\left\{
        Z_{\set{S}, \set{C}} - R_{\set{S}, \set{C}}
        \le
        -\frac{\delta}{2\sigma^2} \cdot n
    \right\}
    \le
    \exp\left\{
        -
        \frac{\delta^2}{16\sigma^4 \cdot R_{\set{S},\set{C}}}
        \cdot n^2
    \right\}.
\]
Since
$R_{\set{S},\set{C}}\le |\set{S}\setminus\set{C}|$
and 
$\delta\ge\frac{8 \sigma^2 \cdot |\set{S} \setminus \set{C}|}{n}$, the choice of $t$ satisfies $t \ge 1$
and the bound is at most
\[
    \exp\left\{
        -\frac{\delta}{2\sigma^2} \cdot n
    \right\},
\]
and hence is dominated by the upper-tail term for $Z_{\set{C}, \set{S}}$.

The bound for the former summand dominates that of the latter summand for all $\delta~\ge~\frac{8 \sigma^2 \cdot |\set{S} \setminus \set{C}|}{n}$. Hence,
\[
    \Pr\left\{
        Z_{\set{C}, \set{S}} - R_{\set{C}, \set{S}} 
        - (Z_{\set{S}, \set{C}} - R_{\set{S}, \set{C}}) 
        > 
        \frac{\delta}{\sigma^2} \cdot n
    \right\}
    \le 
    2 \exp\left\{-\frac{\delta}{8\sigma^2} \cdot n
    \right\}.
\]
Chaining all the presented identities and inequalities gives the upper bound for the probability of interest.

\subsection{Proof of \texorpdfstring{\Cref{lemma:competitor-projected-outcome-deviation}}{}}
\label{appendix:lemma:competitor-projected-outcome-deviation}

Since $\ortprojmat_\set{C}$ is a projection matrix, we have $\ortprojmat_\set{C} X_\set{C} = 0$. Then, by the definition of $Y$ in (\ref{eq:linear-regression}),
\[
    \ortprojmat_\set{C} Y 
    = 
    \ortprojmat_\set{C} X_{\set{S} \setminus \set{C}} \beta_{\set{S} \setminus \set{C}} 
    + \ortprojmat_\set{C} W.
\]
By the definition of the Euclidean norm, we obtain
\begin{equation} 
    \|\ortprojmat_\set{C} Y\|_2^2 - \|\ortprojmat_\set{C} W\|_2^2 
    = 
    \|\ortprojmat_\set{C} X_{\set{S} \setminus \set{C}} \beta_{\set{S} \setminus \set{C}}\|_2^2 
    + 2 W^T \ortprojmat_\set{C} X_{\set{S} \setminus \set{C}} \beta_{\set{S} \setminus \set{C}}.
\end{equation}
Note that 
\[
    \|\ortprojmat_\set{C} X_{\set{S} \setminus \set{C}} \beta_{\set{S} \setminus \set{C}}\|_2^2 
    = 
    \|G_{\set{S}, \set{C}}^{1/2} \beta_{\set{S} \setminus \set{C}}\|_2^2 \cdot n 
    = 
    (1 + c_1) \cdot \delta_{c_1} \cdot n.
\]
So, $\|\ortprojmat_\set{C} Y\|_2^2 - \|\ortprojmat_\set{C} W\|_2^2 \sim \mathcal{N}(\mu, \nu^2)$ where
\begin{equation*}
    \begin{array}{c}
        \mu = (1 + c_1) \cdot \delta_{c_1} \cdot n; \\
        \\
        \nu^2 = 4 \sigma^2 \cdot (1 + c_1) \cdot \delta_{c_1} \cdot n.
    \end{array}
\end{equation*}

Let $E_\set{C} \coloneqq \|\ortprojmat_\set{C} Y\|_2^2 - \|\ortprojmat_\set{C} W\|_2^2$. Then, the probability of interest can be rewritten as
\[
    \Pr\{E_\set{C} \le \delta_{c_1} \cdot n\} 
    = 
    \Pr\left\{
        \frac{E_\set{C} - \mu}{\nu} 
        \le 
        -\sqrt{\frac{c_1^2 \cdot \delta_{c_1}}{4 \sigma^2 \cdot (1 + c_1)} \cdot n}
    \right\}
\]
where $\frac{E_\set{C} - \mu}{\nu} \sim \mathcal{N}(0, 1)$. Then, we can use the classical tail bound shown in \Cref{lemma:gaussian-tail-bound}:
\[
    \Pr\left\{
        \frac{E_\set{C} - \mu}{\nu} 
        \le 
        -c_1 \cdot \sqrt{\frac{\delta_{c_1}}{4 \sigma^2 \cdot (1 + c_1)} \cdot n}
    \right\} 
    \le 
    \exp\left\{
        -\frac{c_1^2 \cdot \delta_{c_1}}{8 \sigma^2 \cdot (1 + c_1)} 
        \cdot n 
    \right\}
\]
which completes the proof.

\section{Proof of \texorpdfstring{\Cref{theorem:probability-of-error}}{}}
\label{appendix:theorem:probability-of-error}

The proof relies on the union bound coupled with the error probability bound established in \hyperlink{corollary:lemma:pairwise-statistic}{the corollary} of \Cref{lemma:pairwise-statistic}. The result is inspired by Theorem 1 in \cite{wainwright2009information} and the proof follows similar strategy as presented in section IV-C of the original paper.

Denote the set of all subsets of size $s$ except $\set{S}$ as $\mathfrak{C} \coloneqq \{\set{C} \colon |\set{C}| = s, \set{C} \ne \set{S}\}$. By the definition of $\widehat{\set{S}}_{n, s}$ in \eqref{eq:best-subset} and the union bound, we have
\[
    \Pr\{\widehat{\set{S}}_{n, s} \ne \set{S}\} 
    \le 
    \Pr\left[
        \cup_{\set{C} \in \mathfrak{C}}
        \{\widehat{\Delta}_{\set{C}, \set{S}} \le 0\} 
    \right] 
    \le 
    \sum\limits_{\set{C} \in \mathfrak{C}}
    \Pr\{\widehat{\Delta}_{\set{C}, \set{S}} \le 0\}.
\]
For brevity, let 
\[
    \psi_n 
    \coloneqq 
    \frac{1}{c} \cdot \frac{\betamin^2}{\sigma^2} \cdot r_{\set{S}, n}.
\]
Bounding each term in the sum on the right hand side through \hyperlink{corollary:lemma:pairwise-statistic}{the corollary} of \Cref{lemma:pairwise-statistic}, we obtain
\[
    \sum\limits_{\set{C} \in \mathfrak{C}}
    \Pr\{\widehat{\Delta}_{\set{C}, \set{S}} \le 0\} 
    \le 
    3 \sum\limits_{\set{C} \in \mathfrak{C}} 
    \exp\left\{
        -\psi_n \cdot |\set{S} \setminus \set{C}| \cdot n
    \right\}.
\]
    
One way to simplify the obtained majorizing sum is to group its elements by values of $|\set{S} \setminus \set{C}|$. Let $N_m \coloneqq |\{\set{C}\in \mathfrak{C}\colon |\set{S} \setminus \set{C}| = m\}|$. Note that $N_m = \binom{s}{s - m} \cdot \binom{d - s}{m} = \binom{s}{m} \cdot \binom{d - s}{m}$ where we let $\binom{a}{b} = 0$ if $b > a$. Hence, denoting $h = \min\{s, d - s\}$, we can rewrite the sum of the exponential terms as 
\[
    \sum\limits_{\set{C} \in \mathfrak{C}} 
    \exp\left\{
        -\psi_n \cdot |\set{S} \setminus \set{C}| \cdot n
    \right\} 
    = 
    \sum\limits_{m = 1}^{h} 
    N_m \cdot \exp\left\{-\psi_n \cdot m \cdot n\right\}.
\]
We can bound each element of the sum above by the largest one, obtaining 
\begin{equation}\label{eq:probability-of-mistake-bound}
    \begin{split}
        \sum\limits_{\set{C} \in \mathfrak{C}}
        \Pr\{
            \widehat{\Delta}_{\set{C}, \set{S}} \le 0
        \} 
        & \le 
        3 \, h 
        \cdot \max\limits_{m \in \{1, \ldots, h\}}\left\{
            N_m \cdot \exp\left\{- \psi_n \cdot m \cdot n\right\}
        \right\} \\
        & = 
        3 \exp\left\{
            \max\limits_{m \in \{1, \ldots, h\}}\{
                \ln h + \ln N_m - \psi_n \cdot m \cdot n
            \}
        \right\}.
    \end{split}
\end{equation}
For the maximum in the exponent to be non-positive it is sufficient to have
\[
    n 
    \ge 
    \frac{1}{\psi_n} 
    \cdot \max\limits_{m \in \{1, \ldots, h\}} \left\{
        (\ln h + \ln N_m) / m
    \right\}.
\]
Notice that $\max\limits_{m \in \{1, \ldots, h\}} (\ln h + \ln N_m) / m \le \ln h + \max\limits_{m \in \{1, \ldots, h\}} \ln N_m / m$. Since $m \in \{1, \ldots, h\}$, we have 
\[
    \ln N_m / m 
    = 
    \frac{1}{m} 
    \cdot \left[
        \ln \binom{s}{m} + \ln \binom{d - s}{m}
    \right] 
    \le 
    \ln s + \ln (d - s). 
\]

Thus, for any $c_1 \in [0, 1)$ and any $n \ge n_0 + c_1 \cdot n$, where
\begin{equation}\label{eq:n_0}
    n_0 \coloneqq \frac{1}{\psi_n} \cdot (\ln s + \ln (d - s) + \ln h),
\end{equation}
we have
\[
    \sum\limits_{\set{C} \in \mathfrak{C}}
    \Pr\{\widehat{\Delta}_{\set{C}, \set{S}} \le 0\} 
    \le 
    3 \exp\left\{-c_1 \cdot \psi_n \cdot n\right\} 
    = 
    3 \exp\left\{
        -\frac{c_1}{c} 
        \cdot \frac{\betamin^2}{\sigma^2} 
        \cdot r_{\set{S}, n} 
        \cdot n
    \right\}
\]
according to (\ref{eq:probability-of-mistake-bound}) and the definition of $\psi_n$. Note that $m$ disappears from the exponent because of the maximum. Rearranging $n \ge n_0 + c_1 \cdot n$ into $n \ge n_0 / (1 - c_1)$, which is allowed since $c_1 < 1$, we get the following sufficient condition:
\begin{equation}\label{eq:n_0_transformed}
    n 
    \ge 
    n_0 / (1 - c_1) 
    = 
    \frac{c}{1 - c_1} 
    \cdot \frac{\sigma^2}{\betamin^2} 
    \cdot \frac{1}{r_{\set{S}, n}} 
    \cdot (\ln s + \ln (d - s) + \ln h).
\end{equation}
It is only left to make sure that the condition in \Cref{lemma:pairwise-statistic} is satisfied. This is achieved by taking the maximum between $\frac{c \sigma^2}{r_{\set{S}, n} \betamin^2}$ and the lower bound in \eqref{eq:n_0_transformed}, and bounding the result above using $1 - c_1 \le 1$.

\section{Latent-factor design illustration}
\label{appendix:latent-factor-example}

We now give a simple example showing how additional structure in the covariate distribution can make the empirical subset-separability parameter easier to control. In this example, active covariates contain both a low-dimensional factor component and an independent idiosyncratic component, whereas inactive covariates are generated only from the low-dimensional factor component. The argument highlights the following mechanism. First, since all inactive covariates lie in the same low-dimensional nuisance subspace, the residual projection associated with any competing subset can be lower-bounded by projecting out a larger space that does not depend on the particular false replacement variables. Second, after this enlargement, the remaining random matrix is generated only by the idiosyncratic components of the active covariates, so a standard smallest singular value bound applies. Third, in the part of the union bound corresponding to competitors that miss $m$ active variables, the usual factor $\binom{s}{m}\binom{d-s}{m}$ is replaced by $\binom{s}{m}$. This reduction is the source of the weaker sample-size requirement for controlling empirical subset separability in this example.

For a technically clean illustration, we impose two non-disruptive assumptions. The regime that we consider is sparse, that is 
\[
    s \ll d.
\]
Moreover, we assume that there exists $c_2 \in (0, 1)$ such that the factor space dimensionality, denoted as $v$, and the support size together satisfy 
\[
    s + v \le c_2 \, n.
\]

Let the true support be fixed as $\set{S} = \{1, \ldots, s\}$. Consider a Gaussian design generated as follows. Let $Z$ and $U$ be respectively $n \times v$ and $n \times s$ random matrices with independent standard Gaussian entries, and $A \in \RR^{v \times s}$ and $B \in \RR^{v \times (d - s)}$ be deterministic loading matrices. Additionally, let $\tau > 0$ be a constant. Define
\begin{equation}\label{eq:factor-model}
    X_\set{S}=Z A + \tau U,
    \qquad
    X_{\set{S}^c} = Z B.
\end{equation}
Thus the active covariates have a low-dimensional factor component and an idiosyncratic component, while all inactive covariates lie in the same factor subspace of rank at most $v$. In particular,
\[
    \col(X_{\set{S}^c}) \subseteq \col(Z),
\]
where $\col(M)$ denotes the column span of matrix $M$.
The rows of $X$ are i.i.d. Gaussian with covariance matrix
\[
    \Sigma
    =
    \begin{pmatrix}
        A^T A + \tau^2 I_s & A^T B\\
        B^T A & B^T B
    \end{pmatrix}.
\]

Let $\set{C}$ be a competing subset with $|\set{C}| = s$ and $\set{C} \ne \set{S}$, and denote
\[
    h \coloneqq \min\{s, d - s\}.
\]
Fix $m \in \{1, \ldots, h\}$ which corresponds to the number of missing active covariates, and thus to the number of falsely added covariates as well. Write
\[
    \set{J} = \set{S} \setminus \set{C},
    \qquad
    \set{T} = \set{S} \cap \set{C},
    \qquad
    \set{F} = \set{C} \setminus \set{S},
\]
so that $|\set{J}| = |\set{F}| = m$ and $\set{C} = \set{T} \cup \set{F}$. Since $X_\set{F}$ is a submatrix of $X_{\set{S}^c} = Z B$, we have
\[
    \col(X_\set{C})
    =
    \col(X_\set{T},X_\set{F})
    \subseteq
    \col(X_\set{T},Z).
\]
Denote the projection operator onto the orthogonal complement to $\col(X_\set{T}, Z)$ as $\ortprojmat_{[X_\set{T}, Z]}$. Then, $\ortprojmat_\set{C} \succeq \ortprojmat_{[X_\set{T}, Z]}$ since the latter operator projects away the larger space $\col(X_\set{T}, Z)$.
Consequently,
\[
    X_\set{J}^T \ortprojmat_\set{C} X_\set{J}
    \succeq
    X_\set{J}^T \ortprojmat_{[X_\set{T}, Z]} X_\set{J}.
\]
Importantly, the right-hand side does not depend on the false replacement subset $\set{F}$.

Using the model definition in \eqref{eq:factor-model},
\[
    X_\set{J} = Z A_\set{J} + \tau U_\set{J}.
\]
Since $Z$ is included in the projected-out space, the factor component is removed from $X_\set{J}$ when projecting it with $\ortprojmat_{[X_\set{T}, Z]}$. Moreover,
\[
    \col(X_\set{T}, Z)
    =
    \col(U_\set{T}, Z),
\]
because $X_\set{T} = Z A_\set{T} + \tau U_\set{T}$, $U_\set{T}$ is independent of $Z$, and $\tau \ne 0$. Hence,
\[
    \ortprojmat_{[X_\set{T},Z]} X_\set{J}
    =
    \ortprojmat_{[U_\set{T}, Z]} (Z A_\set{J} + \tau U_\set{J})
    =
    \tau \ortprojmat_{[U_\set{T}, Z]} U_\set{J}.
\]
It follows that
\[
    \frac{1}{n}
    X_\set{J}^T \ortprojmat_\set{C} X_\set{J}
    \succeq
    \frac{\tau^2}{n}
    U_\set{J}^T
    \ortprojmat_{[U_\set{T}, Z]}
    U_\set{J}.
\]

Now fix a set $\set{J} \subseteq \set{S}$ with $|\set{J}| = m$. The matrix $U_\set{J}$ is independent of $[U_{\set{T}}, Z]$. If $n > s + v$, then
\[
    \operatorname{rank}([U_{\set{T}}, Z])
    =
    s - m + v
\]
with probability one, and therefore the residual projector $\ortprojmat_{[U_\set{T}, Z]}$ has rank
\[
    q_m = n - s + m - v.
\]
By rotational invariance,
\[
    U_\set{J}^T
    \ortprojmat_{[U_{\set{T}}, Z]}
    U_\set{J}
    \stackrel{d}{=}
    G_m^T G_m,
\]
where $G_m \in \RR^{q_m \times m}$ has independent standard Gaussian entries.

The standard lower-tail bound for the smallest singular value of a Gaussian matrix gives, for every $t \ge 0$,
\[
    \Pr\left\{
        \mineig(G_m^T G_m)
        \le
        \left(\sqrt{q_m} - \sqrt{m} - t\right)_+^2
    \right\}
    \le
    e^{-t^2/2}.
\]
Therefore, if we denote
\[
    Q_m
    \coloneqq
    \min_{\set{C} \subseteq \covars\colon |\set{C}| = s,\ |\set{S} \setminus \set{C}| = m}
    \mineig
    \left(
        \frac{1}{n}
        X_{\set{S} \setminus \set{C}}^T
        \ortprojmat_\set{C}
        X_{\set{S} \setminus \set{C}}
    \right),
\]
then for each fixed $m = 1, \ldots, h$ and every $t \ge 0$,
\begin{equation}\label{eq:Q_m-bound}
    \Pr\left\{
        Q_m
        \le
        \frac{\tau^2}{n}
        \left(\sqrt{n - s + m - v} - \sqrt{m} - t\right)_+^2
    \right\}
    \le
    \binom{s}{m} e^{-t^2/2}.
\end{equation}

Importantly, the union bound is only over the missed true subset $\set{J} \subseteq \set{S}$, not over the false replacement subset $\set{F} \subseteq \set{S}^c$. The latter affects the bound only through the rank-$v$ factor space that is projected out.

Applying the union bound to the lower-tail bounds for $Q_m$ over $m = 1, \ldots, h$, we obtain, for arbitrary $t_1, \ldots, t_h \ge 0$,
\begin{equation}\label{eq:r-union-bound}
    \Pr\left\{
        r_{\set{S}, n}
        \le
        \min_{1 \le m \le h}
        \frac{\tau^2}{n}
        \left(\sqrt{n - s + m - v} - \sqrt{m} - t_m\right)_+^2
    \right\}
    \le
    \sum_{m = 1}^h
    \binom{s}{m} e^{-t_m^2/2}.
\end{equation}

Since $m \le h \le s$ and $n > s + v$, the function
$m \mapsto \sqrt{n - s - v + m} - \sqrt m$ is decreasing. Hence,
\[
    \sqrt{n - s - v + m} - \sqrt m
    \ge
    \sqrt{n - v} - \sqrt s.
\]
Note that $\sqrt{n - v} - \sqrt s = \frac{n - s - v}{\sqrt{n - v} + \sqrt s}$. Since we assume $s + v \le c_2 \, n$, we have $n - s - v \ge (1 - c_2) \cdot n$ and $s \le c_2 \, n$. Hence
\[
    \sqrt{n - v} - \sqrt s
    =
    \frac{n - s - v}{\sqrt{n - v} + \sqrt{s}}
    \ge
    \frac{(1 - c_2) \, n}{(1 + \sqrt c_2) \sqrt n}
    =
    (1 - \sqrt c_2) \sqrt n,
\]
and therefore
\[
    \sqrt{n - s - v + m} - \sqrt m
    \ge
    (1 - \sqrt c_2) \sqrt n
\]
for every $m = 1, \ldots, h$.

Bounding the threshold in \eqref{eq:Q_m-bound} below under $s + v \le c_2 \, n$ and choosing, for any $c_3 \in (0,\, (1 - \sqrt c_2)^2)$, $t =  \sqrt{c_3 n}$ gives
\[
    \Pr\left\{
        Q_m
        \le
        \left(1 - \sqrt c_2 - \sqrt{c_3}\right)^2 \cdot \tau^2
    \right\}
    \le
    \binom{s}{m}\exp\left\{-\frac{c_3 n}{2}\right\}.
\]
Substituting this into the union bound in \eqref{eq:r-union-bound}, we obtain
\[
    \Pr\left\{
        r_{\set{S}, n}
        \le
        \left(1 - \sqrt c_2 - \sqrt{c_3}\right)^2 \cdot \tau^2
    \right\}
    \le
    \left(
        \sum_{m=1}^{h} \binom{s}{m}
    \right)
    \cdot \exp\left\{-\frac{c_3 n}{2}\right\}
    \le
    2^s \exp\left\{-\frac{c_3 n}{2}\right\}
    =
    \exp\left\{- \frac{c_3 n}{2} + s \ln 2\right\}.
\]

Thus the empirical subset-separability parameter is bounded away from zero with high probability once $n$ is of order at least $s + v$. Namely, defining
\[
    n_0 \coloneqq
    \max\left\{
        \frac{1}{c_2} \cdot (s + v),
        \frac{2\ln 2}{c_3} \cdot s
    \right\},
\]
we can take $n \ge n_0 + c_1 n$ for any $c_1 \in [0, 1)$, and then
\[
    \Pr\left\{
        r_{\set{S}, n}
        \le
        \left(1 - \sqrt c_2 - \sqrt{c_3}\right)^2 \cdot \tau^2
    \right\}
    \le
    \exp\left\{-c_1 \cdot \frac{c_3 n}{2}\right\}.
\]
Rearranging $n \ge n_0 + c_1 n$ gives us the final sufficient condition:
\[
    n 
    \ge 
    \frac{1}{1 - c_1} 
    \cdot \max\left\{
        \frac{1}{c_2} \cdot (s + v),\ 
        \frac{2\ln 2}{c_3} \cdot s
    \right\}.
\]

This is advantageous whenever the latent rank $v$ is small compared with the full combinatorial complexity of choosing false replacements. In the sparse regime $s \ll d$, we have
\[
    \ln\binom{d - s}{s}
    \asymp
    s \ln\left(\frac{d - s}{s}\right),
\]
so a requirement of order $s + v$ is smaller than the full log-combinatorial price whenever
\[
    v \ll s \ln\left(\frac{d - s}{s}\right).
\]

This calculation illustrates the role of the latent-factor structure. Although there are $d - s$ inactive variables, all inactive covariates lie in the same sample subspace of rank at most $v$.
Hence different choices of the false replacement set $\set{F}\subseteq\set{S}^c$ do not change the enlarged projected-out space $\col(X_\set{T},Z)$ used in the lower bound. The price paid for the inactive covariates is instead the loss of $v$ residual degrees of freedom. By contrast, for an unstructured random design with independent inactive covariates, a union-bound argument would typically range over all choices of $\set{F}$, producing the additional factor $\binom{d - s}{m}$ at distance $m$.

Thus, in this stylized random-design model, the empirical subset-separability parameter can be controlled without paying the full log-combinatorial price associated with all false replacement labels. This provides a concrete example of the phenomenon suggested by the fixed-design theorem: the relevant complexity is not necessarily the number of candidate supports, but the number of distinguishable projection directions generated by the design.

\section{Restricted candidate class illustration}
\label{appendix:restricted-candidate-class}

We now give a second illustration, showing how restrictions on the candidate model class can reduce the price of controlling the empirical subset-separability parameter. The mechanism is different from the latent-factor example in \Cref{appendix:latent-factor-example}. There, the reduction comes from the geometry of the design: many inactive covariates lie in the same low-dimensional nuisance subspace. Here, the candidate supports are restricted to a smaller family. For example, such restrictions may come from oracle or expert information, or from modeling constraints such as grouped variables. As a result, both the empirical separability control and the error probability union bound range over the restricted candidate class, rather than over all subsets of size $s$.

Let the rows of $X$ be independent draws from $\mathcal{N}(0, \Sigma)$, where $\Sigma \in \RR^{d \times  d}$ is positive definite. Let $\set{S} \subseteq \covars$ be the true support, with $|\set{S}| = s$. Suppose that the search space is restricted to a candidate family
\[
    \mathfrak{S}
    \subseteq
    \{\set{C} \subseteq \covars \colon |\set{C}| = s\},
    \qquad
    \set{S} \in \mathfrak{S}.
\]
The corresponding restricted empirical subset-separability parameter is
\[
    r_{\set{S}, n}(\mathfrak{S})
    \coloneqq
    \min_{\set{C} \in \mathfrak{S}\colon \set{C} \ne \set{S}}
    \mineig
    \left(
        \frac{1}{n}
        X_{\set{S} \setminus \set{C}}^T
        \ortprojmat_\set{C}
        X_{\set{S} \setminus \set{C}}
    \right).
\]
This is the analogue of $r_{\set{S}, n}$ when best subset selection is performed only over $\mathfrak{S}$. 

For a given $\set{C}\in\mathfrak{S}$, the population analogue of the residualized Gram matrix is the conditional covariance
\[
    \Sigma_{\set{S} \setminus \set{C} \mid \set{C}}
    \coloneqq
    \Sigma_{\set{S} \setminus \set{C}, \set{S} \setminus \set{C}}
    -
    \Sigma_{\set{S} \setminus \set{C}, \set{C}}
    \Sigma_{\set{C},\set{C}}^{-1}
    \Sigma_{\set{C},\set{S} \setminus \set{C}},
\]
and the restricted population-level subset separability is denoted as
\[
    \rho_{\set{S}}(\mathfrak{S})
    \coloneqq
    \min_{\set{C}\in\mathfrak{S} \colon \set{C} \ne \set{S}}
    \mineig
    \left(
        \Sigma_{\set{S} \setminus \set{C} \mid \set{C}}
    \right).
\]

We first show how pairwise subset separation can be controlled given a fixed competitor $\set{C} \in \mathfrak{S}$ with $m = |\set{S} \setminus \set{C}|$. By the Gaussian conditioning formula,
\[
    X_{\set{S} \setminus \set{C}}^T
    \ortprojmat_\set{C}
    X_{\set{S} \setminus \set{C}}
    \stackrel{d}{=}
    \Sigma_{\set{S} \setminus \set{C} \mid \set{C}}^{1/2}
    G_{\set{C}}^T G_{\set{C}}
    \Sigma_{\set{S} \setminus \set{C} \mid \set{C}}^{1/2}
\]
where $G_{\set{C}}$ is an $(n - s) \times m$ random matrix with independent standard Gaussian entries. For the random matrix on the right-hand side, we have the deterministic bound
\[
    \mineig
    \left(
        \Sigma_{\set{S} \setminus \set{C} \mid \set{C}}^{1/2}
        G_{\set{C}}^T G_{\set{C}}
        \Sigma_{\set{S} \setminus \set{C} \mid \set{C}}^{1/2}
    \right)
    \ge
    \mineig(\Sigma_{\set{S} \setminus \set{C} \mid \set{C}})
    \cdot
    \sigma_{\min}^2(G_{\set{C}})
    \ge
    \rho_{\set{S}}(\mathfrak{S})
    \cdot
    \sigma_{\min}^2(G_{\set{C}}).
\]
Consequently, for every $g \ge 0$,
\[
    \Pr\left\{
        \mineig
        \left(
            \frac{1}{n}
            X_{\set{S} \setminus \set{C}}^T
            \ortprojmat_\set{C}
            X_{\set{S} \setminus \set{C}}
        \right)
        \le
        g
    \right\}
    \le
    \Pr\left\{
        \frac{\rho_{\set{S}}(\mathfrak{S})}{n} 
        \cdot\sigma_{\min}^2(G_{\set{C}})
        \le 
        g
    \right\}.
\]

Using the standard lower-tail bound for the smallest singular value of a Gaussian matrix, for every $t \ge 0$,
\[
    \Pr\left\{
        \sigma_{\min}(G_{\set{C}})
        \le
        \sqrt{n - s} - \sqrt m - t
    \right\}
    \le
    e^{-t^2/2}.
\]
Hence,
\[
    \Pr\left\{
        \mineig
        \left(
            \frac{1}{n}
            X_{\set{S} \setminus \set{C}}^T
            \ortprojmat_\set{C}
            X_{\set{S} \setminus \set{C}}
        \right)
        \le
        \frac{\rho_{\set{S}}(\mathfrak{S})}{n}
        \left(\sqrt{n - s} - \sqrt m - t\right)_+^2
    \right\}
    \le
    e^{-t^2/2}.
\]

By taking the same $t \ge 0$ for each $m$ and applying the union bound, we have 
\[
    \Pr\left\{
        r_{\set{S}, n}(\mathfrak{S})
        \le
        \min_{1 \le m \le h}
        \frac{\rho_{\set{S}}(\mathfrak{S})}{n}
        \left(
            \sqrt{n - s} - \sqrt m - t
        \right)_+^2
    \right\}
    \le
    |\mathfrak{S}| e^{-t^2/2}.
\]
In the unrestricted problem, the corresponding number of competitors is $\binom{d}{s} - 1$.
In the restricted problem, this factor is replaced by $|\mathfrak{S}|$, which may be much smaller.

Suppose, for example, that $n \ge 4s$. Since $m \le h \le s$,
\[
    \sqrt{n - s} - \sqrt{m}
    \ge
    \sqrt{n - s} - \sqrt{s}
    \ge
    \frac{1}{3} \sqrt n.
\]
Thus, for any $c_2 \in [0, 1/3)$, choosing $t = c_2 \sqrt n$ yields
\[
    \Pr\left\{
        r_{\set{S}, n}(\mathfrak{S})
        \le
        \left(\frac{1}{3} - c_2\right)^2
        \rho_{\set{S}}(\mathfrak{S})
    \right\}
    \le
    |\mathfrak{S}|
    \exp\left\{-\frac{c_2^2}{2} \cdot n\right\}
    =
    \exp\left\{
        \ln|\mathfrak{S}| - \frac{c_2^2}{2}\cdot n
    \right\}.
\]
Taking
\[
    n
    \ge
    \frac{1}{1 - c_1}
    \max\left\{
        4s,
        \frac{2}{c_2^2} \ln|\mathfrak{S}|
    \right\}
\]
for any $c_1 \in [0, 1)$, we obtain
\[
    \Pr\left\{
        r_{\set{S}, n}(\mathfrak{S})
        \le
        \left(\frac{1}{3} - c_2\right)^2
        \rho_{\set{S}}(\mathfrak{S})
    \right\}
    \le
    \exp\left\{- \frac{c_1 c_2^2}{2} \cdot n\right\}.
\]

So, for a meaningful random-design result on the basis of \Cref{theorem:probability-of-error}, one would also assume that the candidate family is uniformly separated at the population level, meaning that
\[
    \rho_{\set{S}}(\mathfrak{S})
    > 0.
\]

The simple calculation carried out in this section shows that restricting the candidate class can reduce the price of controlling empirical subset separability. If the search is unrestricted, then the union bound ranges over all possible competing subsets. If the search is restricted to $\mathfrak{S}$, then the relevant price is instead $\ln|\mathfrak{S}|$.

For example, if an initial screening step or expert constraint reduces the admissible covariate set to a subset $\set{M}\subseteq\covars$ with $\set{S}\subseteq\set{M}$ and $|\set{M}|=p$, then
\[
    \mathfrak{S}
    =
    \{\set{C}\subseteq\set{M}\colon |\set{C}|=s\},
    \qquad
    \ln|\mathfrak{S}|
    =
    \ln\binom{p}{s},
\]
instead of $\ln\binom{d}{s}$. As another example, suppose that the covariates are partitioned into $G$ groups of equal size $d_g$, and that admissible supports must be unions of groups. If $s=kd_g$, then selecting $k$ groups gives
\[
    |\mathfrak{S}|=\binom{G}{k},
\]
rather than the unrestricted count $\binom{d}{s}=\binom{Gd_g}{kd_g}$. Thus, the improvement in this illustration comes from restricting the collection of events over which the union bound is applied. The Gaussian design is used only as a standard setting in which the lower tail of the minimum eigenvalue of each relevant residualized Gram matrix can be controlled explicitly.

\section{Auxiliary bounds}

\subsection{Basic tail probability bounds}

\begin{lemma}\label{lemma:gaussian-tail-bound}
    If $X \sim \mathcal{N}(\mu, \nu^2)$, then $\Pr\{X \ge \mu + t\} \le \exp\left\{-\frac{t^2}{2\nu^2}\right\}$ for all $t \ge 0$.
\end{lemma}

The following two bounds appear in \citet{wainwright2009information}, taken from \citet{laurentmassart2000}.

\begin{lemma}\label{lemma:chi-square-upper-tail-bound}
    If $X \sim \chi^2_{k}$, then $\Pr\{X - k \ge 2 \sqrt{k x} + 2x\} \le \exp\{-x\}$.
\end{lemma}

\begin{lemma}\label{lemma:chi-square-lower-tail-bound}
    If $X \sim \chi^2_{k}$, then $\Pr\{X - k \le -2 \sqrt{k x}\} \le \exp\{-x\}$.
\end{lemma}

\subsection{Probability bounds derived from basic tail bounds}

The following two bounds originate from Appendix D of \citet{wainwright2009information}: the former one is directly mentioned, and the latter follows from some derivations. In short, the results follow from \Cref{lemma:chi-square-upper-tail-bound,lemma:chi-square-lower-tail-bound} for the upper and the lower-tail bound respectively, by setting $x = k t$ and using $\sqrt t \le t$ since $t \ge 1$.

\begin{lemma}\label{lemma:chi-square-normalized-upper-tail-bound}
    If $X \sim \chi^2_{k}$, then $\Pr\left\{\frac{X - k}{k} \ge  4t\right\} \le \exp\{-k t\}$ for all $t \ge 1$.
\end{lemma}

\begin{lemma}\label{lemma:chi-square-normalized-lower-tail-bound}
    If $X \sim \chi^2_{k}$, then $\Pr\left\{\frac{X - k}{k} \le  -2 \sqrt{t}\right\} \le \exp\{-k t\}$ for all $t \ge 1$.
\end{lemma}
    
\end{appendices}

\end{document}